\documentclass[a4paper,final,conference,10pt]{IDAACS}%
\usepackage{url}
\usepackage{cite}
\usepackage{amsmath}
\usepackage{amssymb}
\usepackage{graphicx}
\usepackage{amsfonts}%
\usepackage{multirow}
\title{Wavelet Ensemble Kalman Filters}

\begin{document}

\setlength{\oddsidemargin}{0.0in} \setlength{\evensidemargin}{0.0in}
\twocolumn[ \vspace*{0.5in}
\begin{center}{\huge\textbf{Wavelet Ensemble Kalman Filters}}\end{center}\bigskip
\begin{center}Jonathan D. Beezley, Jan Mandel, and Loren Cobb \\University of Colorado Denver, Denver, CO\end{center}\bigskip ]

\let\thefootnote\relax\footnotetext{This research was partially supported by NSF grant
AGS-0835579 and NIH grant RC1-LM010641-01.}

\begin{abstract}
We present a new type of the EnKF for data assimilation in spatial
models that uses diagonal approximation of the state covariance in the wavelet space to achieve
adaptive localization. The efficiency of the new method is demonstrated on an example.
\end{abstract}


\section{Introduction}

Incorporating new data into computations in progress is a well-known problem
in many areas, including weather forecasting, signal processing, and computer
vision. Sequential Bayesian techniques based on the state-space model are
known as filtering or data assimilation
\cite{Anderson-1979-OF,Evensen-2009-DAE}. The probability distribution of the
system state is advanced in time by the computational model, while data is
incorporated from time to time by modifying the probability distribution of
the state by an application the Bayes theorem. Gaussian probability
distributions are represented by their mean and covariance. An assumed
constant state covariance then yields the classical optimal statistical
interpolation (OSI). The Kalman filter (KF) evolves the state covariance, but
it needs to maintain the covariance matrix, so it is not suitable for
high-dimensional systems. The ensemble Kalman filter (EnKF)
\cite{Evensen-2009-DAE} replaces the state covariance by sample covariance of
an ensemble of simulations. The EnKF allows an implementation without any
change to the model; the model only needs to be capable of exporting its state
and restarting from the state modified by the EnKF. Convergence of the
ensemble covariance to the state covariance is guaranteed in the large
ensemble limit \cite{Mandel-2009-CEK,LeGland-2009-LSA}, but a good
approximation may require hundreds of ensemble members \cite{Evensen-2009-DAE}
because the covariance between physically distant variables is small, yet the
sample covariance for a small ensemble has many large long-range terms.
Localization techniques \cite{Furrer-2007-EHP,Hunt-2007-EDA,Anderson-2001-EAK}
improve the approximation by suppressing the long-range terms.

In \cite{Mandel-2010-DDC,Mandel-2010-FFT,Mandel-2010-DAM}, we have proposed an
alternative approach, the fast Fourier transform (FFT) EnKF. The FFT EnKF
assumes that the state is a random field that is approximately homogeneous in
space. Then the covariance matrix in the frequency domain can be well
represented by its diagonal, which gives a good approximation even for very
small ensembles. However, the covariance is not represented well when it
varies with location. The sample covariance can be used for the
cross-covariance between different physical fields in the state
\cite{Mandel-2010-DAM}, which, however, may again cause spurious long-range correlations.

In this paper we extend the spectral approach to the \emph{wavelet EnKF},
resulting in an automatic localization that varies in space adaptively. We
also introduce a new technique for automatic localization of the
cross-covariances, based on a projection and a diagonalization in the spectral
space. The efficiency of the new methods is demonstrated on an example.

\section{Related work}

Diagonal approximation of the covariance in the frequency space was proposed
for weather fields \cite{Berre-2000-ESM}. Wavelets are well suited for
approximation of meteorological fields \cite{Fournier-2000-IOW}, and the
diagonal approximation was extended to wavelet spaces
\cite{Deckmyn-2005-WAR,Fournier-2010-DWM,Pannekoucke-2007-FPW}. The Fourier
domain KF \cite{Castronovo-2008-MTC} is the KF applied to independent
frequency modes. The Laplace operator represented by a diagonal matrix in the
frequency space was used for a fast OSI \cite{Mandel-2010-DAM}. The inverse of
the Laplace operator was proposed as a covariance model
\cite{Kitanidis-1999-GCF}, but higher negative powers \cite{Mandel-2010-DAM}
yield better distributions.

\section{The ensemble Kalman filter}

\label{sec:assimilation}

The modeled quantity is the probability distribution of the state $u$,
represented by an ensemble of simulation states $\left\{  u_{1},\ldots
,u_{N}\right\}  $. A data vector $d$ is linked with the state $u$ by the
observation matrix $H$ such that if the model and the data are correct, then
$d=Hu$. The data error is assumed to have normal distribution with zero mean
and known covariance $R$. When the data arrives, the ensemble is updated by%
\begin{equation}
u_{k}^{a}=u_{k}+Q_{N}H^{\mathrm{T}}\left(  HQ_{N}H^{\mathrm{T}}+R\right)
^{-1}\left(  d+e_{k}-Hu_{k}\right)  ,\label{eq:enkf}%
\end{equation}
where $Q_{N}$ is the sample covariance computed from the ensemble, and the
data perturbation vectors $e_{k}$ are sampled from the data error
distribution. The ensemble members are then advanced in time by the model
until the next data vector is to be assimilated. See \cite{Evensen-2009-DAE}
for details.

\section{Orthogonal wavelet transform}

\begin{figure}[ptb]
\vspace*{-0.2in}
\par
\begin{center}
\hspace*{-0.4in} \includegraphics[width=3.8in]{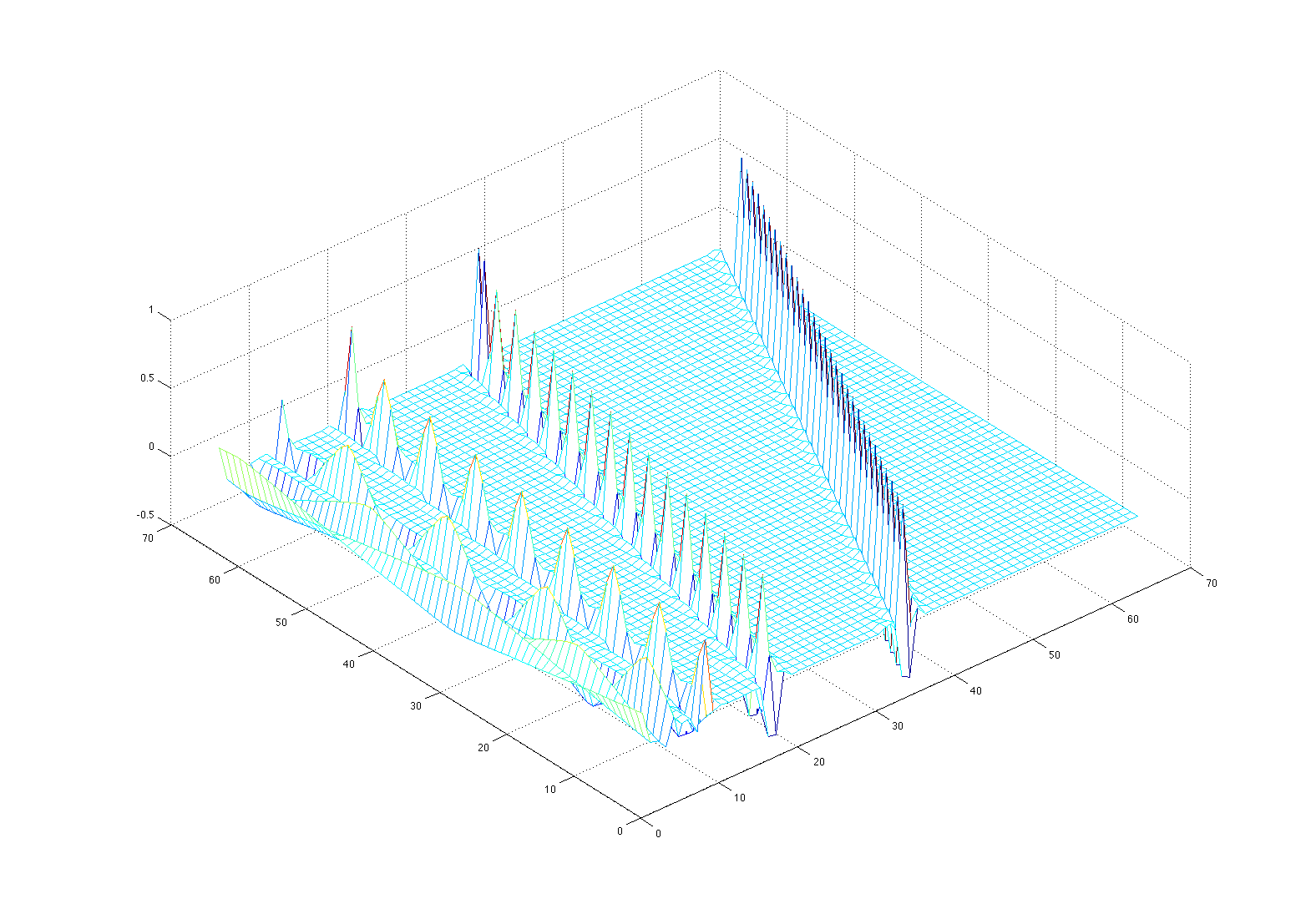}
\end{center}
\caption{Wavelet transform matrix with $n=64$ and $5$ octaves, for the Coiflet
2 wavelet.}%
\label{fig:wavelet-x}%
\end{figure}

For vector $u$, denote $\widehat{u}=Fu$ the transform%
\[
\widehat{u}=\left[  \widehat{u}_{\ell}\right]  _{\ell=1}^{n},\quad\widehat
{u}_{\ell}=f_{\ell}\cdot u,
\]
where $f_{\ell}$ are the rows of $F$ and $F$ is orthonormal, $F^{-1}%
=F^{\mathrm{T}}$. Matrices are then transformed by
\[
\widehat{M}=FMF^{-1}=FMF^{\mathrm{T}}.
\]

In the FFT EnKF, we use the Fourier sine transform. In the wavelet EnKF, we
use the orthogonal wavelet transform, where the rows of $F$ are $f_{\ell
}=f_{\left(  k,j\right)  }$, given by values of the family of functions
$\phi(2^{k}x-j)$ at the points $x=i/n$ for a total $n$ of composite indices
$\ell=\left(  k,j\right)  $ such that $0\leq j\leq2^{k}-1$. Each range of the
indices $\left(  k,j\right)  $ with a fixed $k$ is called an octave. The
function $\phi$ is called the mother wavelet and is chosen so that the rows
$f_{\ell}$ are orthonormal. The indices need to span full octaves, which
restricts the dimension $n$ to a power of $2$ (cf. Fig. \ref{fig:wavelet-x}),
though generalizations are possible. See
\cite{Wickerhauser-1994-AWA,Fournier-2000-IOW} for further details. We use
WaveLab \cite{Buckheit-1995-WRR} to perform the wavelet transform $\widehat
{u}=Fu$ in Matlab. The complexity of the fast wavelet transform is only
$O\left(  n\right)  $, compared to $O\left(  n\log n\right)  $ for the FFT.

In more than 1D, the spectral transformations are applied in each dimension
separately, i.e., the basis functions are taken to be tensor products of 1D
basis functions. Unlike in the case of the FFT, the tensor product of wavelets
creates a certain bias for coordinate directions; the impact of this bias on specific
applications must be examined \cite{Deckmyn-2005-WAR}. However, 2D wavelets
do not seem to be practical yet.

\section{Spectral approximation of the covariance}

\label{sec:spectral-cov}

\subsection{Single variable}

Consider first the case when the model state consists of one variable, in 1D
only. Denote by $u\left(  x_{i}\right)  $, $i=1,\ldots,n$ the entry of the
vector $u$, corresponding to node $x_{i}$. If the random field $u$ is
stationary, then the covariance matrix satisfies $Q\left(  x_{i},x_{j}\right)
=c\left(  x_{i}-x_{j}\right)  $ for some covariance function $c$, and
matrix-vector multiplication $v=Qu$ is the convolution,%
\[
v\left(  x_{i}\right)  =\sum_{j=1}^{n}Q\left(  x_{i},x_{j}\right)  u\left(
x_{j}\right)  =\sum_{j=1}^{n}u\left(  x_{j}\right)  c\left(  x_{i}%
-x_{j}\right)  .
\]
In the spectral domain, convolution becomes entry-by-entry multiplication of
vectors, that is, the multiplication by a diagonal matrix, at least approximately.

Let $\widehat{u}_{ik}$ be the entries of the column vector $\widehat{u}%
_{k}=Fu_{k}$, i.e., ensemble member $k$ in the spectral space. Then we have
the sample covariance in the spectral domain,%
\begin{align}
\widehat{C}(u,u)  &  =FC(u,u)F^{-1}\label{eq:freq-sample-cov}\\
&  =\frac{1}{N-1}{\sum_{k=1}^{N}F}\left(  u_{k}-\overline{u}\right)  \left(
\widehat{u}_{k}-\overline{u}\right)  ^{\mathrm{T}}F^{-1}\nonumber\\
&  =\frac{1}{N-1}{\sum_{k=1}^{N}}\left(  \widehat{u}_{k}-\overline{\widehat
{u}}\right)  \left(  \widehat{u}_{k}-\overline{\widehat{u}}\right)  ,\nonumber
\end{align}
where $\overline{\widehat{u}}=\frac{1}{N}\sum_{k=1}^{N}\widehat{u}_{k}$ is the
sample mean. Assuming that the covariance is primarily a function of the
physical distance between the nodes $x_{i}$, we approximate the forecast
covariance in the spectral domain by the diagonal $\widehat{D}(u,u)$ of
$\widehat{C}(u,u)$,
\begin{align}
\left(  \widehat{D}(u,u)\right)  _{ii}  &  =\frac{1}{N-1}{\sum_{k=1}^{N}%
}\left\vert \widehat{u}_{ik}-\overline{\widehat{u}}_{i}\right\vert
^{2},\label{eq:dii}\\
\left(  \widehat{D}(u,u)\right)  _{ij}  &  =0\text{ for }i\neq j.\nonumber
\end{align}

\subsection{Multiple variables}

When the model state consists of more than one variable, the covariance is
split into blocks of cross-covariances between each variable. The diagonal
blocks of the covariance can be approximated as in (\ref{eq:dii}).
Off-diagonal blocks cannot, in general, be treated the same way because the
meshes over which different variables are defined may not coincide. In this
case, we define an interpolation operator $P_{uv}$ that projects a variable
$u$ to the mesh of variable $v$. While this matrix is rectangular in general,
we require an approximate left inverse, such as the Moore-Penrose
pseudoinverse, $P_{uv}^{\dagger}$, where $P_{uv}^{\dagger}P_{uv}\approx I$.
In practice, more sophisticated and efficient methods with a sparse representation of
the approximate left inverse are possible.  In Section \ref{sec:computation}, 
we consider only cross-covariances between identical grids, leaving
a more detailed examination of these techniques for future research.

By projecting $u$ to the grid of $v$, it is possible to proceed as in the
single variable case and construct an approximate cross-covariance,
\[
C(u,v)\approx P_{uv}^{\dagger}C(P_{uv}u,v)\approx P_{uv}^{\dagger}%
F^{-1}\widehat{C}\left(  P_{uv}u,v\right)  F.
\]

\section{Spectral EnKF}

\label{sec:fftenkf}

\subsection{Single variable}

First assume that the observation function $H=I$, that is, the whole state $u$
is observed. The evaluation of the EnKF formula (\ref{eq:enkf}) in the
frequency domain, with the diagonal spectral approximation $\widehat{D}(u,u)$
of the covariance becomes%
\begin{equation}
\widehat{u}^{a}=\widehat{u}+\widehat{D}(u,u)\left(  \widehat{D}(u,u)+\widehat
{R}\right)  ^{-1}\left(  \widehat{d}+\widehat{e}-\widehat{u}\right)  .
\label{eq:freq-enkf}%
\end{equation}
The analysis ensemble is obtained by the inverse transform at the end,
$u_{k}^{a}=F^{-1}\widehat{u}_{k}^{a}.$ Since $\widehat{D}$ is diagonal,
(\ref{eq:freq-enkf}) can be evaluated very efficiently in several important
cases: \emph{(i)} $\widehat{R}$ is diagonal; then $\widehat{D}+\widehat{R}$ is
also diagonal and the evaluation of (\ref{eq:freq-enkf}) reduces to
term-by-term operations on vectors; \emph{(ii) }$\widehat{R}$ is the perturbed
data sample covariance and we can use the Sherman-Morrison-Woodbury formula, which only
requires solving a system of dimension equal to the ensemble size $N$;
\emph{(iii)} $\widehat{R}$ is approximated by the diagonal part of the sample covariance in the
spectral domain, just like the state covariance.

\subsection{Multiple variables}

\label{sec:multiple}

Consider the state with multiple variables and the covariance and the
observation matrix in the block form. Using the spectral covariance estimation
in Sec. \ref{sec:spectral-cov},%
\begin{equation}
QH^{\mathrm{T}}=\left[  C(u_{i},H_{j}u)\right]  \approx P^{\dagger}%
F^{-1}\left[  \widehat{D}(P_{i}u_{i},H_{j}u)\right]  F,\nonumber
\end{equation}
with $P_{i}$ the interpolation operator from mesh $i$ to the observation grid,
$P^{\dagger}$ the block diagonal matrix made up of $P_{1}^{\dagger}%
,\cdots,P_{n_{v}}^{\dagger}$, and $F$ consisting likewise of the spectral
transform in each variable. Similarly, we have $HQH^{\mathrm{T}}\approx
F^{-1}\left[  \widehat{D}(H_{i}u,H_{j}u)\right]  F$, where the blocks are
diagonal, so the term $HQH^{\mathrm{T}}+R$ in (\ref{eq:enkf}) can be inverted
easily in the spectral domain. Since the cross covariances in the term
$QH^{\mathrm{T}}$ in (\ref{eq:enkf}) are not involved in a matrix inversion,
we can use another approximation, such the usual sample covariance
\cite{Mandel-2010-DAM}, which may be more suitable when the field depends on
the observed field non-locally, such as by advection.


\section{Computational example}
\label{sec:computation}
\begin{figure}[ptb]
\par
\begin{center}
\hspace*{-0.15in}
\begin{tabular}
[c]{cc}%
Variable $1$ & Variable $2$\\
\includegraphics[width=1.5in]{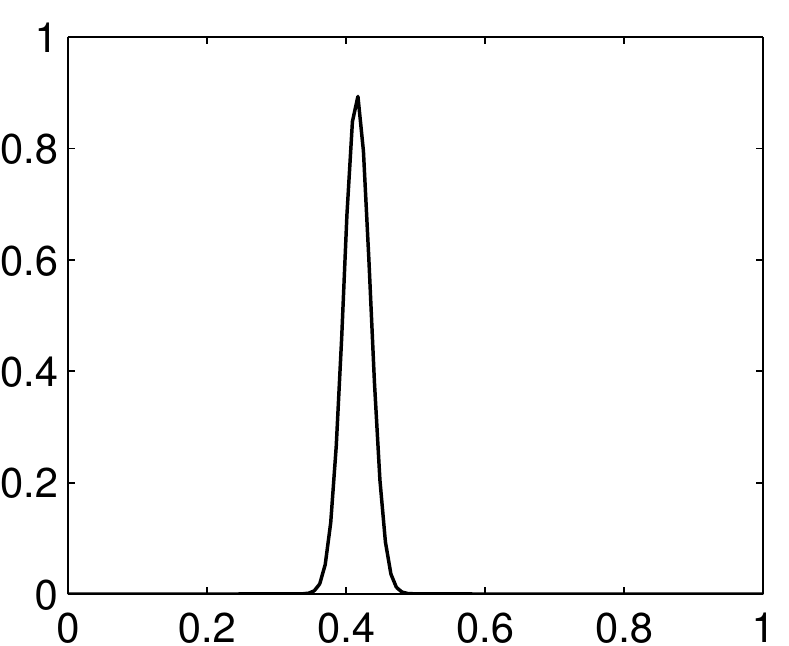} &
\includegraphics[width=1.5in]{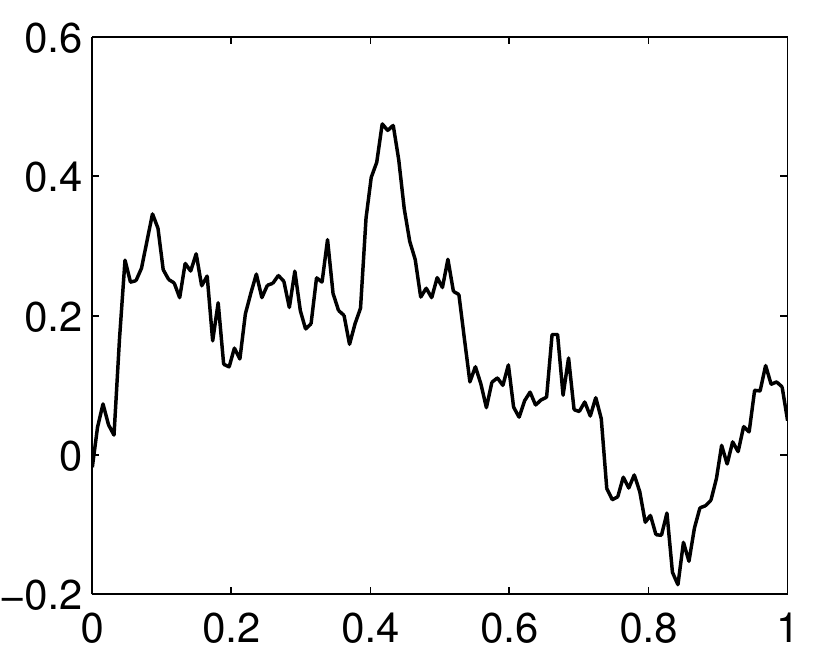}\\
(a) & (b)
\end{tabular}
\end{center}
\caption{A single realization of the stochastic model.}%
\label{fig:model_vars}%
\end{figure}

We conclude with a simple example highlighting the advantages of the proposed
method. We create a synthetic two-variable model state in one spatial
dimension. Both variables are discretized on the same mesh of size $128$ over
the domain $[0,1]$. The first variable is a~simple Gaussian shape with random
center, width, and height defined by
\begin{equation}
u_{1}(x)=h\exp(-\left(  x-c\right)  ^{2}/w^{2}),\label{eq:ex_var1}%
\end{equation}
where $c\sim\mathcal{N}(0.3,0.1^{2})$, $w\sim\mathcal{N}(0.1,0.01^{2})$, and
$h\sim\mathcal{N}(1,0.1^{2})$. The second variable is made up of the sum of
two components. The first is a smooth random field made up of a sum of sine
functions with random amplitudes. The amplitudes are sampled from a Gaussian
distribution with variance decaying as the inverse square of the frequency.
The second component is identical to the first variable scaled by $0.3$. The
variables are chosen in this way so that the diagonal of $\operatorname*{Cov}%
(u_{1},u_{1})$ is non-uniform with a peak at $x=0.3$. Further, it is not
obvious looking at a single realization (Figure \ref{fig:model_vars}) that the
two variables are correlated; however, it is easily verified that
$\operatorname*{Cov}(u_{1},u_{2})=0.3\operatorname*{Cov}(u_{1},u_{1})$.

\begin{figure}[ptb]
\par
\begin{center}
\hspace*{-0.15in}
\begin{tabular}
[c]{cc}%
Large sample covariance & Small sample covariance\\
\includegraphics[width=1.25in]{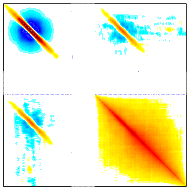} &
\includegraphics[width=1.25in]{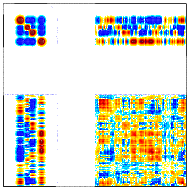}\\
(a) & (b)\\
Covariance by FFT & Covariance by wavelets\\
\includegraphics[width=1.25in]{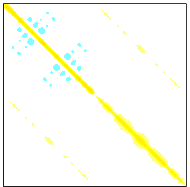} &
\includegraphics[width=1.25in]{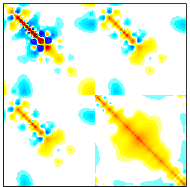}\\
(c) & (d)
\end{tabular}
\end{center}
\caption{In (b)-(d), covariance estimates with an ensemble size $10$ 
are compared with the sample covariance of a large ensemble of size $1000$ in (a).
The color scale is the same in each figure with positive and negative
correlations displayed in red and blue, respectively.}%
\label{fig:model_cov}%
\end{figure}

Using a relatively small ensemble of size $10$, one can see from Figure
\ref{fig:model_cov} that both the FFT\ and wavelet covariance estimates
decay to zero with distance as expected. In comparison to a large sample
covariance (size $1000$), the spectral estimates offer a far better estimate
of the true covariance than the traditional sample covariance when computed
with the same ensemble. In addition, the wavelet covariance estimate correctly
captures the spatial variability of the first variable's distribution. This is
in contrast to the FFT estimate that smooths the peak located at $0.3$
uniformly across the domain.

In order to test the effectiveness of the spectral EnKF itself, we simulate an
observation of the first variable with an observation grid that corresponds to
the discretization of the model variables. We choose an observation generated
from (\ref{eq:ex_var1}) with $c=0.4$, $w=0.12$, and $h=1.5$. We apply the
spectral EnKF described in Section \ref{sec:fftenkf} with spectral
transformations constructed from the discrete sine and Coiflet 2 wavelet
transforms. We use same ensemble of size $10$ from Figure \ref{fig:model_cov}
for each assimilation test and choose a small data covariance $R=0.01^{2}I$ in
order to force the assimilation to react visibly to the observation.

\begin{figure}[ptb]
\par
\begin{center}
\hspace*{-0.1in}
\begin{tabular}
[c]{ccc}%
EnKF & FFT EnKF & Wavelet EnKF\\
\includegraphics[width=1in]{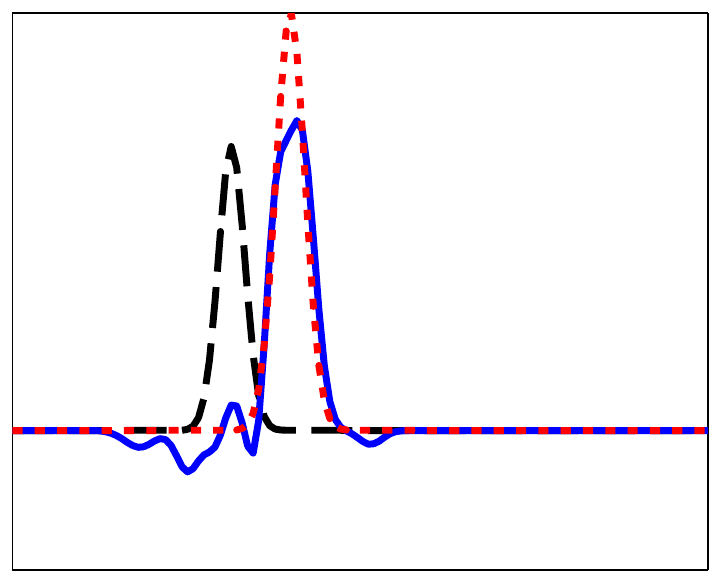} & \hspace{-.2in}
\includegraphics[width=1in]{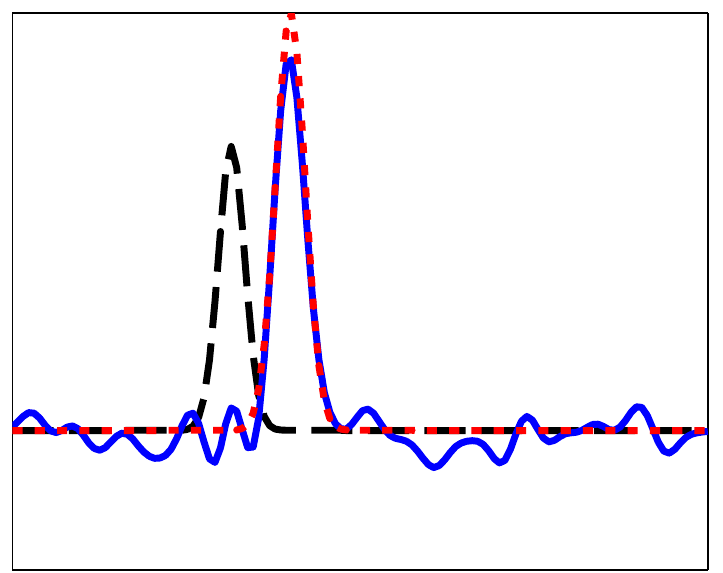} & \hspace{-.2in}
\includegraphics[width=1in]{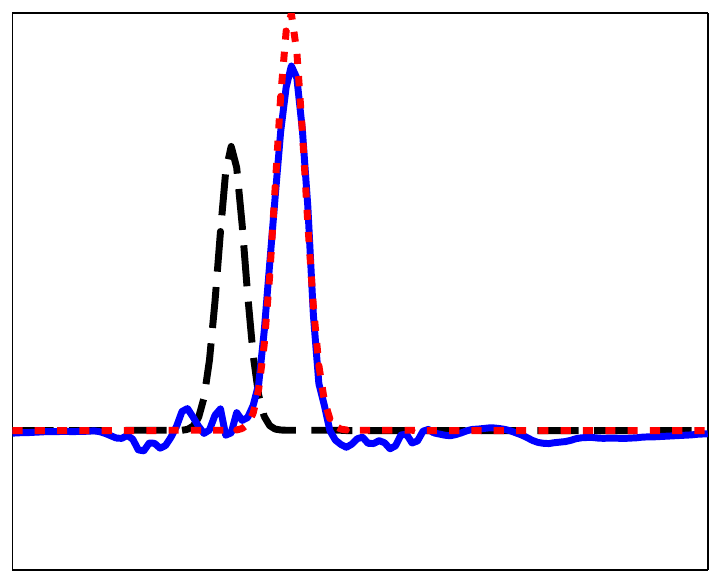}\\
(a) & (b) & (c)\\
\includegraphics[width=1in]{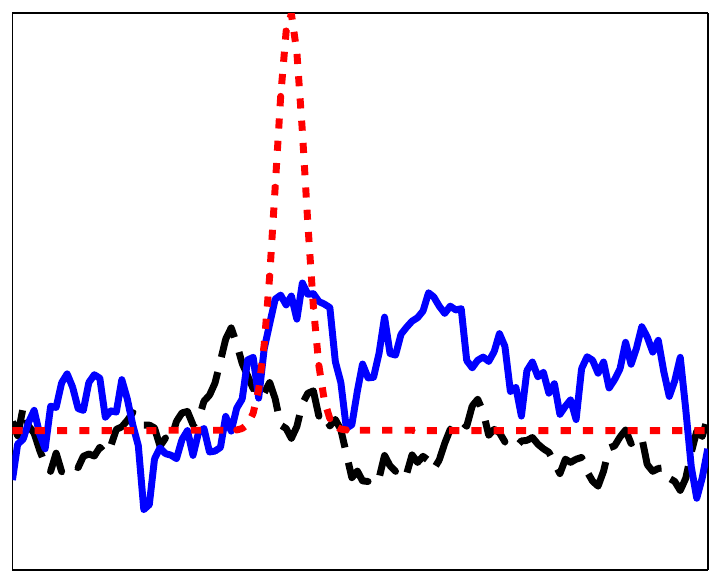} & \hspace{-.2in}
\includegraphics[width=1in]{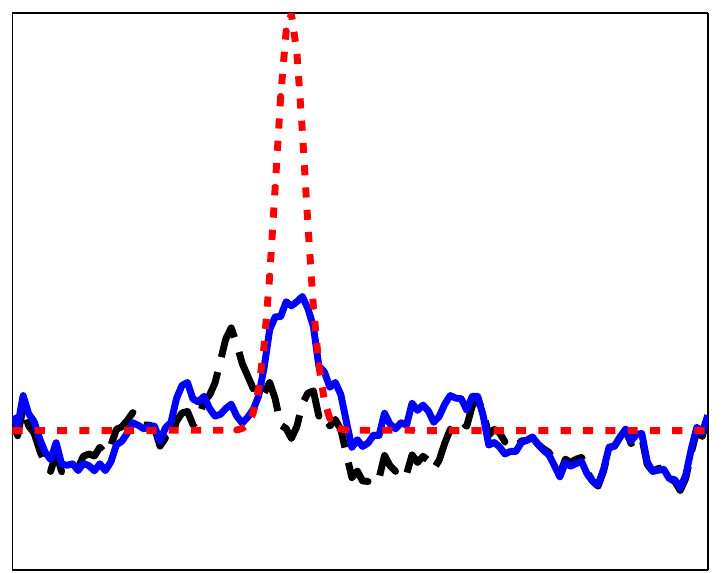} & \hspace{-.2in}
\includegraphics[width=1in]{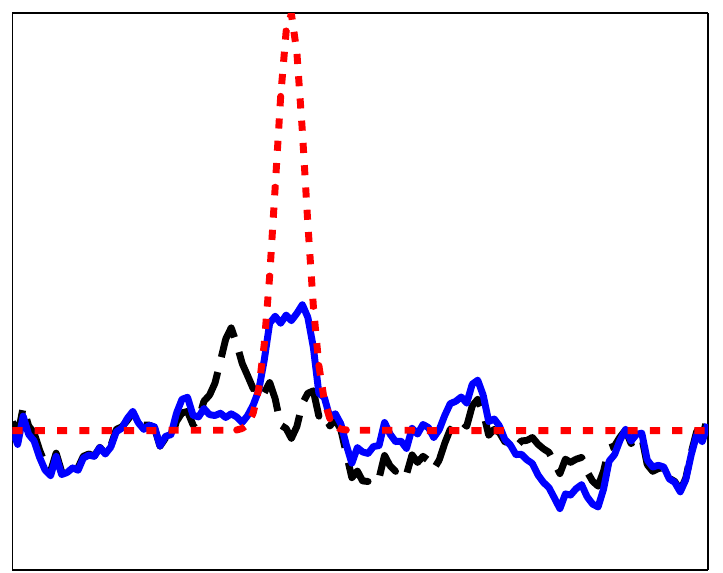}\\
(d) & (e) & (f)
\end{tabular}
\end{center}
\caption{The first (a-c) and second (d-f) variables of a single ensemble
member for each algorithm compared. The analysis is displayed as a solid blue
line along with the forecast as a dashed black line and observation as a
dotted red line.}%
\label{fig:enkf_var1}%
\end{figure}

Figure \ref{fig:enkf_var1} shows the forecast and analysis of the first
ensemble member for each method. The data is displayed in each figure to help
gauge the accuracy of the response of the observed peak at $x=0.4$. The
traditional EnKF produces spurious noise near the peak of the forecast
variable one, while the innovation in variable two is large throughout the
domain rather than local to the observed peak. The FFT EnKF seems to react
well to the observation in the second variable, but the first variable
exhibits spurious noise similar to a Gibb's effect throughout the domain.
Finally, wavelet EnKF appears to provide the best analysis with very little
spurious noise and a seemingly appropriate reaction to the observation.

\section{Conclusion}

Preliminary results indicate that the wavelet EnKF offers a vast improvement
over the traditional assimilation techniques with similar ensemble sizes and
without any localization needed. It requires only spatial information of the
computational mesh, with no expert knowledge necessary to construct background
covariances. Future work will analyze the effects of varying computational
meshes as well as more complex observation functions on a range of operational
computational models.

\bibliographystyle{IEEEtran}
\bibliography{IEEEabrv,../../references/geo,../../references/other,../../references/epi}
\enlargethispage{-1in}  
\end{document}